\documentclass[a4paper,12pt]{article}
\usepackage{latexsym,amsthm}
\newtheorem{define}{Definition}
\newtheorem{lemma}{Lemma}
\newtheorem{theorem}{Theorem}
\newcommand{\Remark}{\par\noindent{\it
Remark.  }}
\def \H {{I\!\!H}}

\begin{document}
\begin{center}
{\Large\bf
Coxeter Decompositions of Bounded Hyperbolic Pyramids and Triangular
Prisms.
 }

\medskip
{\Large A.~Felikson}

\vspace{35pt}
\parbox{10.5cm}
{\scriptsize
{\bf Abstract.}
Coxeter decompositions of hyperbolic simplices where studied in
math.MG/0212010 and math.MG/0210067.
In this paper we use the methods of these works
to classify Coxeter decompositions of
bounded convex pyramids and triangular prisms in the hyperbolic space
$\H^3$.

}
\end{center}
\section*{Introduction}

Let $P$ be a convex polyhedron in the hyperbolic space $H^3$.

\begin{define}
A polyhedron $P$ is called a {\bf Coxeter polyhedron} if all dihedral
angles of $P$ are integer parts of $\pi$.

\end{define}

\begin{define}\label{def1}
A polyhedron $P$ admits a {\bf Coxeter decomposition}
if $P$ can be tiled by finitely many Coxeter polyhedra
such that any two tiles
having a common face are symmetric with respect to this face.

\end{define}

In this paper we classify  Coxeter decompositions of
bounded convex pyramids and triangular prisms in the hyperbolic space
$\H^3$.

\subsection*{Basic definitions}

If $P$  admits a Coxeter decomposition
we also say that $P$ is a {\bf quasi-Coxeter} polyhedron.
The tiles in  Definition~\ref{def1} are called
{\bf fundamental polyhedra} and denoted by $F$.
Clearly, any two fundamental polyhedra are congruent to each other.
A plane $\alpha$ containing a face of a fundamental polyhedron
is called a {\bf mirror} if $\alpha$ contains no face of $P$.

\medskip

In this paper any polyhedron is either a quasi-Coxeter
polyhedron or a polyhedron bounded by the mirrors of some
Coxeter decomposition.
By the "decomposition" we always mean a Coxeter decomposition.

\begin{define}
Given a Coxeter decomposition of a polyhedron $P$,
a {\bf dihedral angle} of $P$ formed up by facets
$\alpha$ and $\beta$ is called {\bf fundamental}
if no mirror contains $\alpha \cap \beta$.
In this case the {\bf edge} $\alpha \cap \beta$ of $P$ is called
{\bf fundamental} too. \\
A {\bf vertex} $A$ is called {\bf fundamental}
if no mirror contains $A$.

\end{define}

\section{Fundamental polyhedron}

Let $P$ be a bounded quasi-Coxeter polyhedron.
From now on $P$ is either
a bounded pyramid or a bounded triangular prism in $\H^3$.

\vspace{6pt}

\noindent
{\bf Notation.}\\
Denote by $\alpha \cap \beta$
the intersection of
the sets $\alpha$ and $\beta$ in the {\it inner}
part of $\H^3$.\\
Let $\alpha$ be a face of a polyhedron, denote by $\overline
\alpha$ a plane containing $ \alpha$.

\begin{lemma}\label{inter}

Let $F$ be a Coxeter polyhedron.\\
Then either $F$ is a tetrahedron or $F$ has two faces $\alpha$
and $\beta$ such that
$\overline \alpha \cap \overline \beta =\emptyset$.

\end{lemma}

\begin{proof}

For any  polyhedron without obtuse
angles we have
$$ \alpha \cap  \beta =\emptyset
\quad \Rightarrow \quad
\overline \alpha \cap \overline \beta =\emptyset$$
(see~\cite{A2}).
Any Coxeter polyhedron has no obtuse angles.
Thus, it is enough
to prove that $F$ has faces
$\alpha$ and $\beta$ such that
$\alpha \cap \beta =\emptyset.$

Suppose that $\alpha \cap \beta \ne \emptyset$
for any faces $\alpha$ and $\beta$ of $F$.
Let $A$ be an arbitrary vertex of $F$. Let
$\alpha_1$,...,$\alpha_k$ be all the faces of $F$
containing $A$.
Let $\beta$ be any face of $F$ such that $A \notin \beta$.
By the assumption,
$\alpha_i \cap \beta \ne \emptyset \quad \forall i=1,...,k$.
Therefore, $F$ has no faces except $\alpha_i$ and $\beta$,
i.e. $F$ is a pyramid.

Suppose that $A$ is not an ideal vertex of $F$.
Then $A$ is incident exactly to three faces of $F$
and $F$ is a tetrahedron.
Suppose that $A$ is an ideal vertex of $F$ and
 the pyramid $F$ is not a tetrahedron.
Consider two faces (say, $\alpha_1$ and $\alpha_3$) having no
common edges.
Since $A\notin \H^3$, we have
$\alpha_1 \cap \alpha_3=\emptyset$.

\end{proof}

\begin{lemma}
\label{l.3.1}
If $P$ is a bounded pyramid, then $F$ is a tetrahedron.

\end{lemma}

\begin{proof}
Suppose that $F$ is not a tetrahedron.
By Lemma~\ref{inter}, $F$ has two faces $\alpha$ and $ \beta$
such that $\overline \alpha \cap \overline \beta =\emptyset$.
Let $F_0$ be a fundamental polyhedron  in $P$;
let $\alpha_0$ and $\beta_0$ be its disjoint faces.
Consider a sequence of fundamental polyhedra $F_i\in P$, $i\in
\mathbf Z$, such that
$\alpha_i=\alpha_{i+1}$, if $i$ is odd, and
$\beta_i=\beta_{i+1}$, if $i$ is even
(see Fig.~\ref{seq1}).

The sequence is finite,
since $P$ contains finitely many fundamental polyhedra.
The polyhedra $F_i$ cannot make a cycle,
since $\alpha_i\cap \beta_i =\emptyset$.
Let $F_k$ and $F_s$ be the endpoints
of the sequence.  Then $\alpha_k$ or $\beta_k$ belongs to
some face $\gamma$  of $P$.
Similarly, $\alpha_s$ or $\beta_s$ belongs to
some face $\delta$ of $P$.
Obviously, $\overline \alpha_i \cap \overline \beta_j=\emptyset$
for any
$i,j$. But $\gamma$ intersects $\delta$, since $P$ is a bounded
pyramid.  The contradiction shows that $F$ is a tetrahedron.

\begin{figure}[!h]
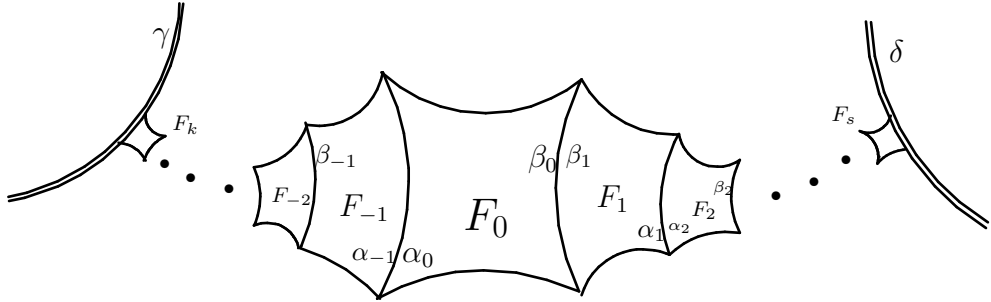

\begin{center}
\input pic/seq1.tex
\end{center}
\caption{A sequence of the fundamental polyhedra.}
\label{seq1}
\end{figure}

\end{proof}

\noindent
{\bf Remark.} The idea of the proof of Lemma~\ref{l.3.1}
belongs to O.~V.~Schwarzman.

\begin{figure}[!h]
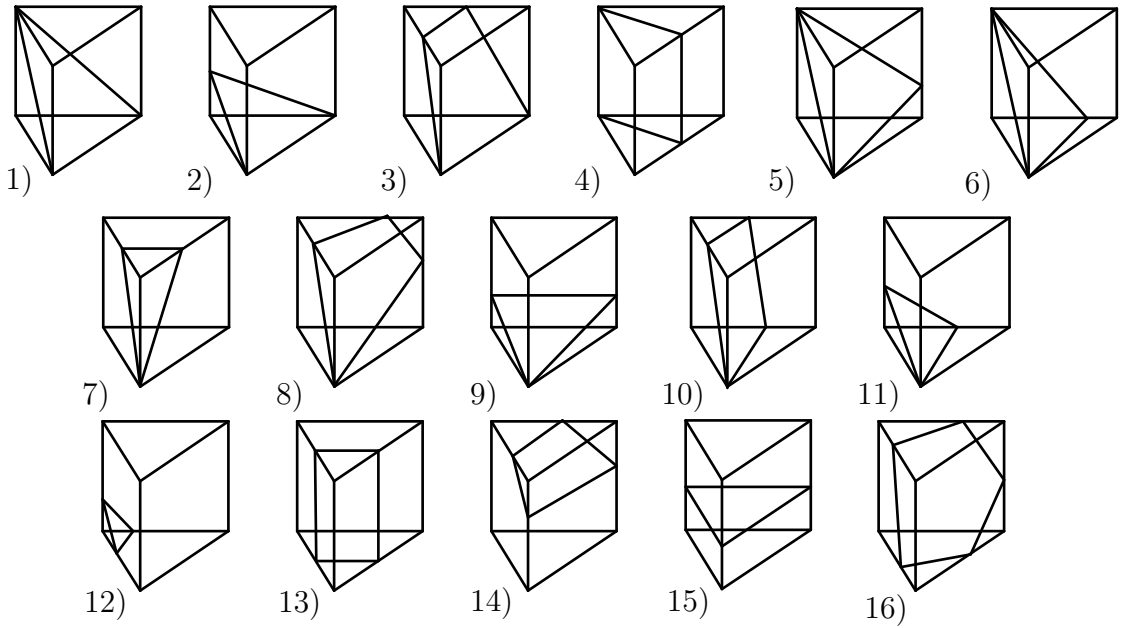

\begin{center}
\input pic/mirrors.tex
\end{center}
\caption{Mirrors in the prism.}
\label{mirrors}
\end{figure}

\vspace{7pt}
Let $\alpha$ be a mirror in a Coxeter decomposition of the
triangular prism $P$. Then
$\alpha$  intersects $P$ as  shown in one
of  Fig.~\ref{mirrors}.1--\ref{mirrors}.16.

\begin{define}\label{mirr}

We say that a mirror is {\bf pentagonal} if it intersects the prism as
shown in Fig.~\ref{mirrors}.16.

\end{define}

\begin{define}
We say that a triangular quasi-Coxeter prism $P$ is {\bf
minimal} if any prism $P'$ inside $P$
is fundamental.

\end{define}

\begin{lemma}
\label{l5}
Let $P$ be a minimal prism and $F$ be a fundamental polyhedron.
Suppose that $F$ is neither a tetrahedron nor a triangular prism.
Then any mirror is pentagonal and any dihedral angle of $P$ is
fundamental.

\end{lemma}

\begin{proof}

Suppose that the decomposition contains a tetrahedron or
a pyramid.  Then by Lemma~\ref{l.3.1} $F$ is a tetrahedron.
Hence, $P$ contains neither pyramid nor tetrahedron.
Since $P$ is minimal, $P$ contains no smaller  triangular prisms.
Therefore, any mirror is pentagonal
(see Fig.~\ref{mirrors}).
Clearly, all dihedral angles of $P$ are fundamental.

\end{proof}

\begin{lemma}
Let $P$ be a triangular prism. Then $F$ is either a tetrahedron or
a triangular prism.

\end{lemma}

\begin{proof}
It is sufficient to prove the lemma for a minimal prism $P$.
Suppose that $F$ is neither a tetrahedron nor a triangular prism.
Then by Lemma~\ref{l5} any mirror is pentagonal and any dihedral
angle of $P$ is fundamental.

Let $ABCDE$ be any pentagonal mirror (see Fig.~\ref{5}).
This mirror separates the points
 $M_1$, $N_1$, $N_3$  from the points $M_2$, $M_3$, $N_2$.

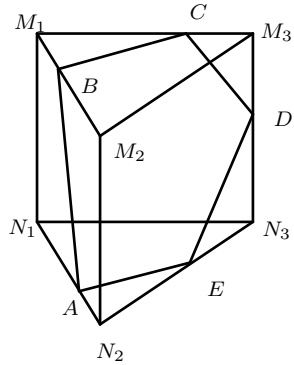
\begin{figure}[!h]
\begin{center}
{\scriptsize
\setlength{\unitlength}{0.075in}
\begin{picture}(18.55,23.68)
\put(0.41,23.00){{\setbox0=\hbox{$M_1$}\lower\ht0\box0}}
\put(17.61,22.33){{\setbox0=\hbox{$M_3$}\lower\ht0\box0}}
\put(7.35,14.00){{\setbox0=\hbox{$M_2$}\lower\ht0\box0}}
\put(6.15,0.00){{\setbox0=\hbox{$N_2$}\lower\ht0\box0}}
\put(17.75,8.50){{\setbox0=\hbox{$N_3$}\lower\ht0\box0}}
\put(0.01,8.66){{\setbox0=\hbox{$N_1$}\lower\ht0\box0}}
\put(13.88,4.33){{\setbox0=\hbox{$E$}\lower\ht0\box0}}
\put(18.55,16.16){{\setbox0=\hbox{$D$}\lower\ht0\box0}}
\put(12.68,23.66){{\setbox0=\hbox{$C$}\lower\ht0\box0}}
\put(5.08,18.50){{\setbox0=\hbox{$B$}\lower\ht0\box0}}
\put(3.75,3.00){{\setbox0=\hbox{$A$}\lower\ht0\box0}}
\special{em:linewidth 0.014in}
\put(3.48,19.16){\special{em:moveto}}
\put(4.95,3.66){\special{em:lineto}}
\put(12.68,5.66){\special{em:lineto}}
\put(17.08,16.00){\special{em:lineto}}
\put(12.41,21.66){\special{em:lineto}}
\put(3.48,19.16){\special{em:lineto}}
\put(17.08,21.66){\special{em:moveto}}
\put(17.08,8.50){\special{em:lineto}}
\put(6.41,14.50){\special{em:moveto}}
\put(6.41,1.33){\special{em:lineto}}
\put(2.01,21.66){\special{em:moveto}}
\put(2.01,8.50){\special{em:lineto}}
\put(2.01,8.50){\special{em:moveto}}
\put(17.08,8.50){\special{em:lineto}}
\put(6.41,1.33){\special{em:lineto}}
\put(2.01,8.50){\special{em:lineto}}
\put(2.01,21.66){\special{em:moveto}}
\put(17.08,21.66){\special{em:lineto}}
\put(6.41,14.50){\special{em:lineto}}
\put(2.01,21.66){\special{em:lineto}}
\end{picture}
}
\end{center}
\caption{A pentagonal mirror.}
\label{5}
\end{figure}

We say that a polyhedron $W$ is {\bf good} if
it contains the points $M_2$, $M_3$, $N_2$
and if it can be cut out of $P$ by a single pentagonal
mirror.
Let $W$ be a minimal good polyhedron
(i.e. no polyhedron inside $W$ is good).
The minimal polyhedron exists, since the decomposition contains
finitely many mirrors. Suppose that $W$ has vertices
$M_2M_3N_2ABCDE$ (see Fig.~\ref{5}).

Suppose that $W$ has a non-fundamental dihedral angle $\alpha$.
By Lemma~\ref{l5} $\alpha$ is an angle formed up by $ABCDE$ and some
face of $P$. Let $\Pi$ be a mirror decomposing the angle
$\alpha$. Clearly, $\Pi$ does not separate the points $M_1$,
$N_1$ and $N_3$ one from another. Since $\Pi$ is pentagonal,
it separates all the points above from the points
$M_2$, $M_3$, $N_2$. Therefore, $\Pi$ cuts out of $P$ a good
polyhedron contained in $W$. This contradicts to the minimal
property of $W$.

Thus, any dihedral angle of $W$ is fundamental,
and $W$ is a Coxeter polyhedron.
The
prolongations of two disjoint faces cannot intersect each other
(see~\cite{A2}).
This property is broken in the polyhedron $W$
(consider the faces $AEN_2$ and $CDM_3$ which prolongations have
a common line $N_1N_3$). The contradiction shows that $F$ is
either a tetrahedron or a triangular prism.

\end{proof}

\section{Decompositions of pyramids.}

This section describes the way to classify all Coxeter
decompositions of bounded hyperbolic pyramids.

Suppose that a pyramid $P$ has only five vertices
$OA_1A_2A_3A_4$ (where $A_1A_2A_3A_4$ is a base of $P$ and $O$ is an
apex). By Lemma~\ref{l.3.1} $F$ is a tetrahedron. Suppose that any edge
$OA_i$ is fundamental. Consider a small sphere $s$ centered in $O$. The
intersection $s\bigcap P$ is a spherical quadrilateral $q$. Any angle of
 $q$ is fundamental. This contradicts to the fact that the sum of the
angles of a spherical quadrilateral should be greater than $2\pi$.

Hence, we can assume that there is a mirror $m$ through $OA_1$.
This mirror decomposes the pyramid into two smaller pyramids.
One of the smaller pyramids is a tetrahedron and another is
either a tetrahedron or a small quadrilateral pyramid.
Thus, any minimal quadrilateral pyramid
consists of two tetrahedra (not necessary fundamental). The
Coxeter decompositions of hyperbolic tetrahedra
are classified in~\cite{arxiv-tetr}.
Therefore, we can find all
decompositions of minimal quadrilateral pyramids.
Using the decompositions of minimal pyramids we
can find the decompositions of greater pyramids too.
So, using this algorithm, it is possible to find all
decompositions of quadrilateral pyramids.

Analogously, it is possible to classify the decompositions of
pyramids with  bigger number of vertices.
Any pyramid $OA_1...A_n$ has a decomposed edge $OA_i$.
So, the pyramid is decomposed into two smaller pyramids which
possible decompositions we should enumerate at the
previous steps.

This procedure realized by a computer leads to a large list of
quadrilateral quasi-Coxeter pyramids, several
pentagonal quasi-Coxeter pyramids, and exactly one hexagonal
quasi-Coxeter pyramid.
See Table~2 for the result.

\section{ Decompositions of prisms into prisms}

In this section both $P$ and $F$ are triangular prisms.
For any prism we say that a triangular face is a
{\bf base} and a quadrilateral face is a {\bf side}.

\begin{lemma}\label{bok}
No base of $F$ belongs to a side of $P$.

\end{lemma}

\begin{proof}
Suppose that a side $s$ of $P$ contains a base of a fundamental
prism $F_1$. Consider a sequence of fundamental prisms
$F_1$, $F_2$,...,$F_n$,...,
such that $F_i\in P$, and $F_i$ and $F_{i+1}$ have a common
base. This chain is finite, since $P$ contains finitely many
fundamental prisms. The prolongations of the bases of $F_1$
have no common points, since $F_1$ is a Coxeter prism.
So, the prolongations of the bases of two different prisms $F_i$
have no common points. Therefore, the prisms $F_1$,
$F_2$,...,$F_n$, cannot make a cycle and the sequence has an
endpoint $F_n$.  One of the bases of $F_n$ belongs to some face
$b$ of $P$ (otherwise the sequence has a prolongation).  The
face $s$ cannot intersect the face $b$, since these faces are
prolongations of the bases of $F_i$.  This contradicts to the
fact that the side $s$ intersects each face of $P$.

\end{proof}

\begin{lemma}  \label{osn}
No side of $F$ belongs to a base of $P$.

\end{lemma}

\begin{proof}
Denote by $b$ and $s$ the bases of $P$.
Suppose that a side $a_0$ of a fundamental prism $F_0$
belongs to $b$. Consider
a sequence of fundamental prisms
$F_{-k}$,...,$F_0$,...$F_l$,  $F_i\in P$,
such that $F_i$ and $F_{i+1}$ have a common base.
The sequence is finite, and the prisms $F_i$ cannot make a cycle.
Let $a_{-k}$ and $b_l$ be the bases of the endpoints
$F_{-k}$ and $F_l$. By Lemma~\ref{bok} the bases
$a_{-k}$ and $b_l$ belong to two different bases of $P$.
If $a_{-k}$ belongs to the base $s$ then the prolongations
of $a_{-k}$ and $a_0$ have a common point
(see Fig.~\ref{chain}). This is impossible,
since the prolongations of the bases have no common points.

\begin{figure}[!h]
\begin{center}
{\scriptsize
\setlength{\unitlength}{0.075in}
\begin{picture}(23.43,22.31)
\put(9.66,19.05){{\setbox0=\hbox{$F_{-l}$}\lower\ht0\box0}}
\special{em:linewidth 0.014in}
\put(8.96,21.45){\special{em:moveto}}
\put(8.96,19.23){\special{em:lineto}}
\put(7.16,20.23){\special{em:moveto}}
\put(7.16,18.03){\special{em:lineto}}
\put(6.41,21.45){\special{em:moveto}}
\put(6.41,19.23){\special{em:lineto}}
\put(6.41,19.23){\special{em:moveto}}
\put(8.96,19.23){\special{em:lineto}}
\put(7.16,18.03){\special{em:lineto}}
\put(6.41,19.23){\special{em:lineto}}
\put(6.41,21.45){\special{em:moveto}}
\put(8.96,21.45){\special{em:lineto}}
\put(7.16,20.23){\special{em:lineto}}
\put(6.41,21.45){\special{em:lineto}}
\put(17.93,11.21){{\setbox0=\hbox{$F_{-l}$}\lower\ht0\box0}}
\put(4.73,10.38){{\setbox0=\hbox{$F_0$}\lower\ht0\box0}}
\put(18.20,21.55){{\setbox0=\hbox{b}\lower\ht0\box0}}
\put(9.13,2.71){{\setbox0=\hbox{s}\lower\ht0\box0}}
\put(5.26,6.71){\special{em:moveto}}
\put(5.20,6.75){\special{em:lineto}}
\put(4.98,6.88){\special{em:moveto}}
\put(4.91,6.91){\special{em:lineto}}
\put(4.71,7.05){\special{em:moveto}}
\put(4.65,7.08){\special{em:lineto}}
\put(4.43,7.21){\special{em:moveto}}
\put(4.41,7.23){\special{em:lineto}}
\put(4.36,7.26){\special{em:lineto}}
\put(4.16,7.40){\special{em:moveto}}
\put(4.10,7.45){\special{em:lineto}}
\put(3.90,7.60){\special{em:moveto}}
\put(3.83,7.65){\special{em:lineto}}
\put(3.63,7.81){\special{em:moveto}}
\put(3.58,7.86){\special{em:lineto}}
\put(3.40,8.05){\special{em:moveto}}
\put(3.35,8.10){\special{em:lineto}}
\put(3.18,8.28){\special{em:moveto}}
\put(3.13,8.35){\special{em:lineto}}
\put(3.00,8.55){\special{em:moveto}}
\put(2.96,8.61){\special{em:lineto}}
\put(2.83,8.83){\special{em:moveto}}
\put(2.80,8.91){\special{em:lineto}}
\put(2.71,9.16){\special{em:moveto}}
\put(2.68,9.23){\special{em:lineto}}
\put(2.68,9.25){\special{em:lineto}}
\put(2.63,9.50){\special{em:moveto}}
\put(2.63,9.58){\special{em:lineto}}
\put(2.61,9.83){\special{em:moveto}}
\put(2.61,9.91){\special{em:lineto}}
\put(2.65,10.16){\special{em:moveto}}
\put(2.66,10.23){\special{em:lineto}}
\put(2.66,10.25){\special{em:lineto}}
\put(2.71,10.48){\special{em:moveto}}
\put(2.73,10.55){\special{em:lineto}}
\put(2.83,10.76){\special{em:moveto}}
\put(2.85,10.85){\special{em:lineto}}
\put(2.96,11.08){\special{em:moveto}}
\put(3.00,11.15){\special{em:lineto}}
\put(3.13,11.36){\special{em:moveto}}
\put(3.16,11.43){\special{em:lineto}}
\put(3.31,11.63){\special{em:moveto}}
\put(3.36,11.71){\special{em:lineto}}
\put(3.51,11.91){\special{em:moveto}}
\put(3.55,11.96){\special{em:lineto}}
\put(3.56,11.98){\special{em:lineto}}
\put(3.71,12.18){\special{em:moveto}}
\put(3.73,12.20){\special{em:lineto}}
\put(3.76,12.25){\special{em:lineto}}
\put(3.93,12.43){\special{em:moveto}}
\put(3.95,12.45){\special{em:lineto}}
\put(3.98,12.50){\special{em:lineto}}
\put(4.15,12.68){\special{em:moveto}}
\put(4.16,12.70){\special{em:lineto}}
\put(4.20,12.73){\special{em:lineto}}
\put(4.36,12.91){\special{em:moveto}}
\put(4.38,12.93){\special{em:lineto}}
\put(4.41,12.98){\special{em:lineto}}
\put(4.60,13.16){\special{em:moveto}}
\put(4.61,13.18){\special{em:lineto}}
\put(4.66,13.21){\special{em:lineto}}
\put(4.85,13.40){\special{em:moveto}}
\put(4.86,13.41){\special{em:lineto}}
\put(4.90,13.46){\special{em:lineto}}
\put(5.08,13.65){\special{em:moveto}}
\put(5.10,13.66){\special{em:lineto}}
\put(5.15,13.70){\special{em:lineto}}
\put(5.33,13.88){\special{em:moveto}}
\put(5.35,13.90){\special{em:lineto}}
\put(5.40,13.93){\special{em:lineto}}
\put(5.58,14.11){\special{em:moveto}}
\put(5.60,14.13){\special{em:lineto}}
\put(5.63,14.16){\special{em:lineto}}
\put(5.81,14.35){\special{em:moveto}}
\put(5.83,14.36){\special{em:lineto}}
\put(5.86,14.40){\special{em:lineto}}
\put(6.05,14.58){\special{em:moveto}}
\put(6.06,14.60){\special{em:lineto}}
\put(6.10,14.63){\special{em:lineto}}
\put(6.28,14.81){\special{em:moveto}}
\put(6.30,14.83){\special{em:lineto}}
\put(6.33,14.86){\special{em:lineto}}
\put(6.51,15.05){\special{em:moveto}}
\put(6.53,15.06){\special{em:lineto}}
\put(6.56,15.10){\special{em:lineto}}
\put(6.73,15.28){\special{em:moveto}}
\put(6.75,15.30){\special{em:lineto}}
\put(6.78,15.35){\special{em:lineto}}
\put(6.95,15.53){\special{em:moveto}}
\put(7.00,15.58){\special{em:lineto}}
\put(7.15,15.76){\special{em:moveto}}
\put(7.20,15.81){\special{em:lineto}}
\put(7.35,16.01){\special{em:moveto}}
\put(7.40,16.08){\special{em:lineto}}
\put(7.55,16.28){\special{em:moveto}}
\put(7.60,16.35){\special{em:lineto}}
\put(7.73,16.56){\special{em:moveto}}
\put(7.78,16.65){\special{em:lineto}}
\put(7.88,16.88){\special{em:moveto}}
\put(7.90,16.95){\special{em:lineto}}
\put(8.00,17.20){\special{em:moveto}}
\put(8.01,17.28){\special{em:lineto}}
\put(8.06,17.53){\special{em:moveto}}
\put(8.06,17.61){\special{em:lineto}}
\put(8.08,17.86){\special{em:moveto}}
\put(8.06,17.93){\special{em:lineto}}
\put(8.06,17.95){\special{em:lineto}}
\put(8.00,18.18){\special{em:moveto}}
\put(7.96,18.26){\special{em:lineto}}
\put(7.85,18.48){\special{em:moveto}}
\put(7.80,18.55){\special{em:lineto}}
\put(8.33,6.71){\special{em:moveto}}
\put(8.38,6.76){\special{em:lineto}}
\put(8.56,6.93){\special{em:moveto}}
\put(8.61,7.00){\special{em:lineto}}
\put(8.78,7.20){\special{em:moveto}}
\put(8.83,7.25){\special{em:lineto}}
\put(8.98,7.43){\special{em:moveto}}
\put(9.03,7.50){\special{em:lineto}}
\put(9.18,7.70){\special{em:moveto}}
\put(9.21,7.75){\special{em:lineto}}
\put(9.23,7.76){\special{em:lineto}}
\put(9.38,7.96){\special{em:moveto}}
\put(9.41,8.03){\special{em:lineto}}
\put(9.55,8.23){\special{em:moveto}}
\put(9.60,8.30){\special{em:lineto}}
\put(9.73,8.50){\special{em:moveto}}
\put(9.78,8.58){\special{em:lineto}}
\put(9.91,8.80){\special{em:moveto}}
\put(9.96,8.86){\special{em:lineto}}
\put(10.10,9.08){\special{em:moveto}}
\put(10.13,9.15){\special{em:lineto}}
\put(10.26,9.35){\special{em:moveto}}
\put(10.30,9.41){\special{em:lineto}}
\put(10.43,9.63){\special{em:moveto}}
\put(10.46,9.68){\special{em:lineto}}
\put(10.48,9.70){\special{em:lineto}}
\put(10.61,9.91){\special{em:moveto}}
\put(10.65,9.98){\special{em:lineto}}
\put(10.78,10.20){\special{em:moveto}}
\put(10.81,10.26){\special{em:lineto}}
\put(10.95,10.46){\special{em:moveto}}
\put(11.00,10.53){\special{em:lineto}}
\put(11.15,10.73){\special{em:moveto}}
\put(11.18,10.80){\special{em:lineto}}
\put(11.33,11.00){\special{em:moveto}}
\put(11.38,11.06){\special{em:lineto}}
\put(11.53,11.28){\special{em:moveto}}
\put(11.55,11.30){\special{em:lineto}}
\put(11.58,11.35){\special{em:lineto}}
\put(11.75,11.53){\special{em:moveto}}
\put(11.81,11.60){\special{em:lineto}}
\put(12.00,11.78){\special{em:moveto}}
\put(12.05,11.83){\special{em:lineto}}
\put(12.25,12.00){\special{em:moveto}}
\put(12.31,12.05){\special{em:lineto}}
\put(12.53,12.20){\special{em:moveto}}
\put(12.60,12.23){\special{em:lineto}}
\put(12.83,12.33){\special{em:moveto}}
\put(12.91,12.35){\special{em:lineto}}
\put(13.16,12.40){\special{em:moveto}}
\put(13.25,12.40){\special{em:lineto}}
\put(13.50,12.36){\special{em:moveto}}
\put(13.56,12.33){\special{em:lineto}}
\put(13.58,12.33){\special{em:lineto}}
\put(13.81,12.21){\special{em:moveto}}
\put(13.88,12.18){\special{em:lineto}}
\put(14.08,12.03){\special{em:moveto}}
\put(14.15,11.98){\special{em:lineto}}
\put(14.33,11.81){\special{em:moveto}}
\put(14.38,11.76){\special{em:lineto}}
\put(14.55,11.58){\special{em:moveto}}
\put(14.60,11.51){\special{em:lineto}}
\put(14.75,11.31){\special{em:moveto}}
\put(14.78,11.25){\special{em:lineto}}
\put(14.91,11.05){\special{em:moveto}}
\put(14.96,10.98){\special{em:lineto}}
\put(15.10,10.76){\special{em:moveto}}
\put(15.13,10.70){\special{em:lineto}}
\put(15.26,10.50){\special{em:moveto}}
\put(15.30,10.43){\special{em:lineto}}
\put(15.41,10.21){\special{em:moveto}}
\put(15.43,10.18){\special{em:lineto}}
\put(15.45,10.15){\special{em:lineto}}
\put(15.55,9.93){\special{em:moveto}}
\put(15.58,9.85){\special{em:lineto}}
\put(15.68,9.61){\special{em:moveto}}
\put(15.71,9.55){\special{em:lineto}}
\put(4.88,8.41){\special{em:moveto}}
\put(7.80,8.35){\special{em:lineto}}
\put(6.43,6.00){\special{em:moveto}}
\put(9.35,5.95){\special{em:lineto}}
\put(4.81,5.05){\special{em:moveto}}
\put(7.73,4.98){\special{em:lineto}}
\put(7.73,4.98){\special{em:moveto}}
\put(7.80,8.35){\special{em:lineto}}
\put(9.35,5.95){\special{em:lineto}}
\put(7.73,4.98){\special{em:lineto}}
\put(4.81,5.05){\special{em:moveto}}
\put(4.88,8.41){\special{em:lineto}}
\put(6.43,6.00){\special{em:lineto}}
\put(4.81,5.05){\special{em:lineto}}
\put(17.80,9.81){\special{em:moveto}}
\put(17.80,6.86){\special{em:lineto}}
\put(15.40,8.20){\special{em:moveto}}
\put(15.40,5.26){\special{em:lineto}}
\put(14.40,9.81){\special{em:moveto}}
\put(14.40,6.86){\special{em:lineto}}
\put(14.40,6.86){\special{em:moveto}}
\put(17.80,6.86){\special{em:lineto}}
\put(15.40,5.26){\special{em:lineto}}
\put(14.40,6.86){\special{em:lineto}}
\put(14.40,9.81){\special{em:moveto}}
\put(17.80,9.81){\special{em:lineto}}
\put(15.40,8.20){\special{em:lineto}}
\put(14.40,9.81){\special{em:lineto}}
\put(23.43,22.31){\special{em:moveto}}
\put(23.43,7.86){\special{em:lineto}}
\put(6.86,14.40){\special{em:moveto}}
\put(6.86,0.00){\special{em:lineto}}
\put(0.00,22.31){\special{em:moveto}}
\put(0.00,7.86){\special{em:lineto}}
\put(0.00,7.86){\special{em:moveto}}
\put(23.43,7.86){\special{em:lineto}}
\put(6.86,0.00){\special{em:lineto}}
\put(0.00,7.86){\special{em:lineto}}
\put(0.00,22.31){\special{em:moveto}}
\put(23.43,22.31){\special{em:lineto}}
\put(6.86,14.40){\special{em:lineto}}
\put(0.00,22.31){\special{em:lineto}}
\end{picture}
}
\end{center}
\caption{No side of $F$ belongs to a base of $P$.}
\label{chain}
\end{figure}
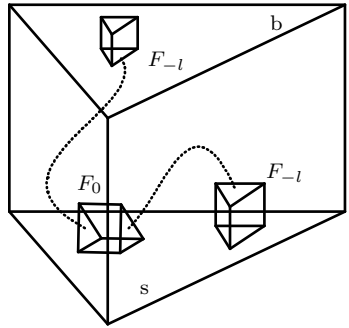

\end{proof}

\begin{lemma}\label{krishki}
Let $b$ be a base of $P$,  $b_1$ and $b_2$ be the bases of
$F$.  Suppose that $b$ is non-fundamental.  Then the tiling of
$b$ by the faces of $F$ is a Coxeter decomposition. A
fundamental polyhedron of this decomposition is a triangle
$b_1$ or $b_2$.

\end{lemma}

\begin{proof}
By Lemma~\ref{osn} any tile is a triangle.
The adjacent tiles of the tiling are the faces of $F$
having a common edge. Thus, if one of the tiles is
a base $b_1$ then the other tiles are the copies of this
base. Clearly, this tiling is a Coxeter decomposition.

\end{proof}

By Lemma~\ref{bok} any side $s$ of $P$ is tiled by the sides of $F$.
This tiling is not necessary a Coxeter decomposition,
but Lemma~\ref{lat}~and~\ref{l_boka} show that this tiling
is very simple.

\begin{figure}[!h]
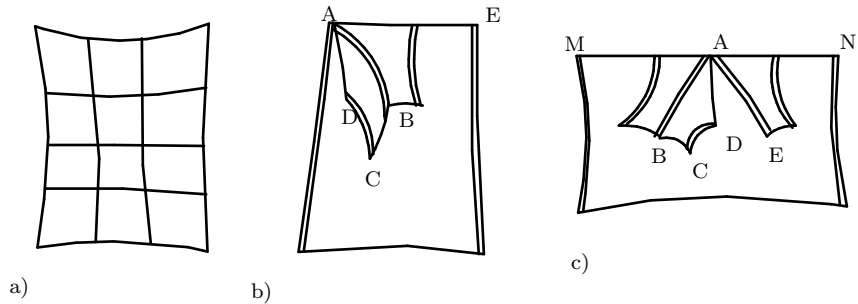

\begin{center}
{\scriptsize
\input pic/latticed.tex
}
\end{center}
\caption{A tiling of a side looks like a lattice.}
\label{lattice}
\end{figure}

\begin{lemma}
\label{lat}
Let $s$ be a side of $P$.
Then the tiling of $s$ is a "lattice"
(see Fig.~\ref{lattice}a), i.e. the tiling
has the following properties:\\
\phantom{iii}
(1) any vertex of the face $s$ belongs to a unique tile;\\
\phantom{iii}
(2) any point of any edge of the quadrilateral $s$
belongs to either one or two tiles.\\
\phantom{iii}
(3) any point in the inner part of $s$ belongs to one, two or four tiles.

\end{lemma}

\begin{proof}

We say that an edge $AB$ of a tile is {\bf horizontal }
if $AB$ belongs to a base of a fundamental prism;
otherwise we say that $AB$ is {\bf vertical}.
Let $ABCD$ be any tile which vertical edges are $AB$ and $CD$.
Consider a sequence of tiles such that
this sequence contains $ABCD$, and
any two neighboring tiles have a common horizontal edge.
We say that this sequence is {\bf vertical}.
By the same way construct a {\bf horizontal} sequence
where the neighboring tiles have a common vertical edge.
Evidently, the vertical sequence ends at two different horizontal
edges of $s$ and the horizontal sequence ends at two different
vertical edges. Let us show that this condition is broken if
the lemma is false.

1). Let $A$ be a vertex of the face $s$. Suppose that $A$
belongs to more then one tile (see Fig.~\ref{lattice}b).
Consider a vertical sequence through the tile $ABCD$.
No pair of tiles from this sequence can be separated
from each other by the line $AB$.
Therefore, this sequence cannot reach the horizontal edge $AE$,
and (1) is proved.

2). Let $A$ be an inner point of a horizontal edge $MN$.
Suppose that $A$ belongs to more then two tiles.
Then one of these tiles has no edges in $MN$
(see Fig.~\ref{lattice}c).
At least two lines through $A$ are vertical edges of some tiles.
Consider a vertical sequence through the tile $ABCD$.
It is evident that the angle $\angle BAE$ contains any tile
from this sequence. Thus, the sequence cannot reach the horizontal
edge $MN$.

The same reasoning works for the points on the
vertical edges, and (2) is proved.

3). Let $A$ be an inner point of the tile $s$. Clearly, $A$
cannot be incident to exactly three tiles. Suppose that $A$ is a vertex of
five or more tiles.

Let $a_1$ be any tile and  $a_1$,...,$a_k$
be an upper part of the vertical sequence containing $a_1$
(i.e. $a_k$ is the upper endpoint
of the sequence). We say that $a_1$ is a tile of the level $k$.
We say that a point $A$ is a point of the level $k$, if $A$ belongs
to a tile of the level $k$. If $A$ belongs to several tiles,
we define the level of $A$ as a maximum.

The part 2) of this proof establishes the lemma for the points
of level 0.
The same reasoning
shows that if the lemma is true for the tiles of
the level  $k$ then it is true for the tiles of the level $k+1$.
Therefore, the lemma is true for a tile of any level.

\end{proof}

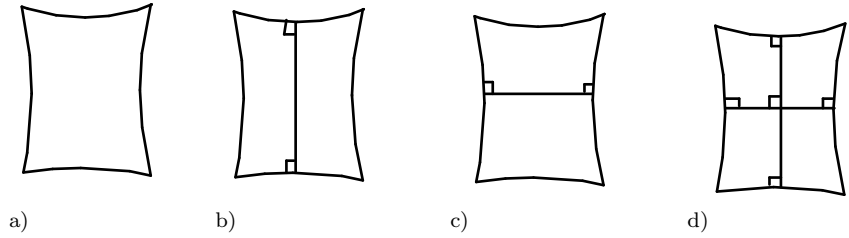
\begin{figure}[!h]
\begin{center}
{\scriptsize
\setlength{\unitlength}{0.075in}
\begin{picture}(58.36,14.41)
\put(47.35,0.00){{\setbox0=\hbox{d)}\lower\ht0\box0}}
\put(30.81,0.00){{\setbox0=\hbox{c)}\lower\ht0\box0}}
\put(14.41,0.00){{\setbox0=\hbox{b)}\lower\ht0\box0}}
\put(0.01,0.00){{\setbox0=\hbox{a)}\lower\ht0\box0}}
\special{em:linewidth 0.014in}
\put(53.88,2.33){\special{em:moveto}}
\put(53.08,2.33){\special{em:lineto}}
\put(53.08,1.83){\special{em:lineto}}
\put(53.21,12.16){\special{em:moveto}}
\put(53.21,11.50){\special{em:lineto}}
\put(53.88,11.50){\special{em:lineto}}
\put(57.48,7.83){\special{em:moveto}}
\put(56.81,7.83){\special{em:lineto}}
\put(56.81,7.16){\special{em:lineto}}
\put(50.01,7.83){\special{em:moveto}}
\put(50.95,7.83){\special{em:lineto}}
\put(50.95,7.33){\special{em:lineto}}
\put(53.08,7.16){\special{em:moveto}}
\put(53.08,8.00){\special{em:lineto}}
\put(53.75,8.00){\special{em:lineto}}
\put(40.68,8.83){\special{em:moveto}}
\put(40.15,8.83){\special{em:lineto}}
\put(40.15,8.33){\special{em:lineto}}
\put(33.21,9.00){\special{em:moveto}}
\put(33.75,9.00){\special{em:lineto}}
\put(33.75,8.16){\special{em:lineto}}
\put(19.35,13.33){\special{em:moveto}}
\put(19.21,12.33){\special{em:lineto}}
\put(19.88,12.33){\special{em:lineto}}
\put(19.35,2.66){\special{em:moveto}}
\put(19.35,3.50){\special{em:lineto}}
\put(19.88,3.50){\special{em:lineto}}
\put(50.01,7.16){\special{em:moveto}}
\put(57.61,7.16){\special{em:lineto}}
\put(53.88,12.16){\special{em:moveto}}
\put(53.88,1.66){\special{em:lineto}}
\put(33.21,8.16){\special{em:moveto}}
\put(40.68,8.16){\special{em:lineto}}
\put(20.01,13.16){\special{em:moveto}}
\put(20.01,2.66){\special{em:lineto}}
\put(58.30,1.13){\special{em:moveto}}
\put(57.68,4.50){\special{em:lineto}}
\put(57.53,7.10){\special{em:lineto}}
\put(57.78,10.48){\special{em:lineto}}
\put(58.36,13.08){\special{em:lineto}}
\put(49.41,1.35){\special{em:moveto}}
\put(53.36,1.66){\special{em:lineto}}
\put(56.76,1.43){\special{em:lineto}}
\put(58.30,1.13){\special{em:lineto}}
\put(49.30,12.93){\special{em:moveto}}
\put(49.85,9.51){\special{em:lineto}}
\put(49.98,6.86){\special{em:lineto}}
\put(49.80,3.75){\special{em:lineto}}
\put(49.41,1.35){\special{em:lineto}}
\put(58.36,13.08){\special{em:moveto}}
\put(56.86,12.56){\special{em:lineto}}
\put(55.18,12.26){\special{em:lineto}}
\put(53.70,12.18){\special{em:lineto}}
\put(51.71,12.33){\special{em:lineto}}
\put(49.70,12.80){\special{em:lineto}}
\put(49.30,12.93){\special{em:lineto}}
\put(41.50,1.80){\special{em:moveto}}
\put(40.98,4.33){\special{em:lineto}}
\put(40.73,7.70){\special{em:lineto}}
\put(40.73,7.76){\special{em:lineto}}
\put(40.86,10.31){\special{em:lineto}}
\put(41.50,13.71){\special{em:lineto}}
\put(41.56,13.75){\special{em:lineto}}
\put(32.61,2.01){\special{em:moveto}}
\put(34.65,2.26){\special{em:lineto}}
\put(36.60,2.33){\special{em:lineto}}
\put(40.03,2.10){\special{em:lineto}}
\put(41.50,1.80){\special{em:lineto}}
\put(32.50,13.60){\special{em:moveto}}
\put(33.05,10.18){\special{em:lineto}}
\put(33.18,7.53){\special{em:lineto}}
\put(32.86,3.38){\special{em:lineto}}
\put(32.61,2.01){\special{em:lineto}}
\put(41.56,13.75){\special{em:moveto}}
\put(40.53,13.35){\special{em:lineto}}
\put(38.58,12.95){\special{em:lineto}}
\put(36.91,12.85){\special{em:lineto}}
\put(34.90,13.00){\special{em:lineto}}
\put(32.90,13.46){\special{em:lineto}}
\put(32.50,13.60){\special{em:lineto}}
\put(24.70,2.13){\special{em:moveto}}
\put(24.18,4.66){\special{em:lineto}}
\put(23.93,8.03){\special{em:lineto}}
\put(23.93,8.10){\special{em:lineto}}
\put(24.06,10.65){\special{em:lineto}}
\put(24.70,14.05){\special{em:lineto}}
\put(24.76,14.08){\special{em:lineto}}
\put(15.81,2.35){\special{em:moveto}}
\put(17.85,2.60){\special{em:lineto}}
\put(19.80,2.66){\special{em:lineto}}
\put(23.23,2.43){\special{em:lineto}}
\put(24.70,2.13){\special{em:lineto}}
\put(15.70,13.93){\special{em:moveto}}
\put(16.25,10.51){\special{em:lineto}}
\put(16.38,7.86){\special{em:lineto}}
\put(16.06,3.71){\special{em:lineto}}
\put(15.81,2.35){\special{em:lineto}}
\put(24.76,14.08){\special{em:moveto}}
\put(23.73,13.68){\special{em:lineto}}
\put(21.78,13.28){\special{em:lineto}}
\put(20.11,13.18){\special{em:lineto}}
\put(18.10,13.33){\special{em:lineto}}
\put(16.10,13.80){\special{em:lineto}}
\put(15.70,13.93){\special{em:lineto}}
\put(9.90,2.46){\special{em:moveto}}
\put(9.28,5.83){\special{em:lineto}}
\put(9.13,8.43){\special{em:lineto}}
\put(9.26,10.98){\special{em:lineto}}
\put(9.90,14.38){\special{em:lineto}}
\put(9.96,14.41){\special{em:lineto}}
\put(1.01,2.68){\special{em:moveto}}
\put(3.05,2.93){\special{em:lineto}}
\put(5.00,3.00){\special{em:lineto}}
\put(8.43,2.76){\special{em:lineto}}
\put(9.90,2.46){\special{em:lineto}}
\put(0.90,14.26){\special{em:moveto}}
\put(1.45,10.85){\special{em:lineto}}
\put(1.58,8.20){\special{em:lineto}}
\put(1.26,4.05){\special{em:lineto}}
\put(1.01,2.68){\special{em:lineto}}
\put(9.96,14.41){\special{em:moveto}}
\put(8.93,14.01){\special{em:lineto}}
\put(6.98,13.61){\special{em:lineto}}
\put(5.31,13.51){\special{em:lineto}}
\put(3.30,13.66){\special{em:lineto}}
\put(1.30,14.13){\special{em:lineto}}
\put(0.90,14.26){\special{em:lineto}}
\end{picture}
}
\end{center}
\caption{Tilings of sides.}
\label{boka}
\end{figure}

\begin{lemma} \label{l_boka}
Any side of $P$ is tiled as shown in
Fig.~\ref{boka}.

\end{lemma}

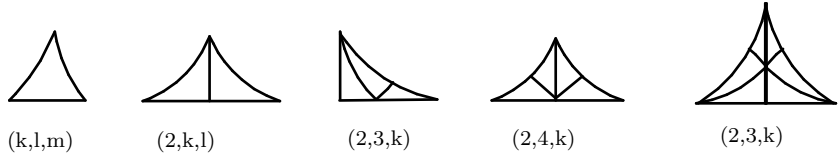
\begin{figure}[!h]
\begin{center}
{\scriptsize
\setlength{\unitlength}{0.075in}
\begin{picture}(57.85,8.93)
\put(0.01,0.00){{\setbox0=\hbox{(k,l,m)}\lower\ht0\box0}}
\special{em:linewidth 0.014in}
\put(5.48,2.16){\special{em:moveto}}
\put(0.15,2.16){\special{em:lineto}}
\put(5.48,2.16){\special{em:moveto}}
\put(4.98,2.81){\special{em:lineto}}
\put(4.30,3.96){\special{em:lineto}}
\put(3.81,5.11){\special{em:lineto}}
\put(3.48,6.26){\special{em:lineto}}
\put(3.35,6.96){\special{em:lineto}}
\put(3.35,6.96){\special{em:moveto}}
\put(2.33,4.93){\special{em:lineto}}
\put(1.56,3.78){\special{em:lineto}}
\put(0.60,2.63){\special{em:lineto}}
\put(0.15,2.16){\special{em:lineto}}
\put(49.75,0.33){{\setbox0=\hbox{(2,3,k)}\lower\ht0\box0}}
\put(35.21,0.16){{\setbox0=\hbox{(2,4,k)}\lower\ht0\box0}}
\put(23.75,0.16){{\setbox0=\hbox{(2,3,k)}\lower\ht0\box0}}
\put(10.41,0.00){{\setbox0=\hbox{(2,k,l)}\lower\ht0\box0}}
\put(53.00,4.51){\special{em:moveto}}
\put(52.23,5.38){\special{em:lineto}}
\put(51.78,5.75){\special{em:lineto}}
\put(53.00,4.51){\special{em:moveto}}
\put(52.08,3.68){\special{em:lineto}}
\put(51.18,3.05){\special{em:lineto}}
\put(50.23,2.56){\special{em:lineto}}
\put(49.20,2.20){\special{em:lineto}}
\put(48.10,1.96){\special{em:lineto}}
\put(53.00,8.93){\special{em:moveto}}
\put(52.70,7.68){\special{em:lineto}}
\put(52.21,6.43){\special{em:lineto}}
\put(51.53,5.18){\special{em:lineto}}
\put(50.58,3.93){\special{em:lineto}}
\put(49.63,2.98){\special{em:lineto}}
\put(48.25,1.96){\special{em:lineto}}
\put(53.00,1.96){\special{em:moveto}}
\put(48.10,1.96){\special{em:lineto}}
\put(53.00,1.96){\special{em:moveto}}
\put(53.00,8.93){\special{em:lineto}}
\put(52.93,4.51){\special{em:moveto}}
\put(53.70,5.38){\special{em:lineto}}
\put(54.16,5.75){\special{em:lineto}}
\put(57.85,1.96){\special{em:moveto}}
\put(56.60,2.25){\special{em:lineto}}
\put(55.38,2.73){\special{em:lineto}}
\put(54.31,3.36){\special{em:lineto}}
\put(53.51,4.00){\special{em:lineto}}
\put(52.93,4.51){\special{em:lineto}}
\put(57.68,1.96){\special{em:moveto}}
\put(56.65,2.68){\special{em:lineto}}
\put(55.60,3.63){\special{em:lineto}}
\put(54.80,4.58){\special{em:lineto}}
\put(54.00,5.85){\special{em:lineto}}
\put(53.43,7.11){\special{em:lineto}}
\put(53.13,8.06){\special{em:lineto}}
\put(52.93,8.93){\special{em:lineto}}
\put(52.93,1.96){\special{em:moveto}}
\put(57.85,1.96){\special{em:lineto}}
\put(52.93,1.96){\special{em:moveto}}
\put(52.93,8.93){\special{em:lineto}}
\put(43.01,2.15){\special{em:moveto}}
\put(42.11,2.50){\special{em:lineto}}
\put(41.05,3.11){\special{em:lineto}}
\put(40.26,3.73){\special{em:lineto}}
\put(39.48,4.55){\special{em:lineto}}
\put(38.90,5.36){\special{em:lineto}}
\put(38.55,5.98){\special{em:lineto}}
\put(38.30,6.48){\special{em:lineto}}
\put(38.30,6.48){\special{em:moveto}}
\put(37.93,5.68){\special{em:lineto}}
\put(37.45,4.88){\special{em:lineto}}
\put(36.80,4.08){\special{em:lineto}}
\put(36.16,3.48){\special{em:lineto}}
\put(35.35,2.88){\special{em:lineto}}
\put(34.53,2.43){\special{em:lineto}}
\put(33.83,2.15){\special{em:lineto}}
\put(33.83,2.15){\special{em:moveto}}
\put(43.01,2.15){\special{em:lineto}}
\put(36.56,3.85){\special{em:moveto}}
\put(38.30,2.31){\special{em:lineto}}
\put(40.16,3.85){\special{em:lineto}}
\put(38.30,6.33){\special{em:moveto}}
\put(38.30,2.15){\special{em:lineto}}
\put(25.76,2.23){\special{em:moveto}}
\put(26.58,2.98){\special{em:lineto}}
\put(26.95,3.43){\special{em:lineto}}
\put(25.76,2.23){\special{em:moveto}}
\put(25.13,2.88){\special{em:lineto}}
\put(24.50,3.75){\special{em:lineto}}
\put(23.90,4.90){\special{em:lineto}}
\put(23.56,5.76){\special{em:lineto}}
\put(23.33,6.63){\special{em:lineto}}
\put(23.26,6.96){\special{em:lineto}}
\put(30.03,2.23){\special{em:moveto}}
\put(28.61,2.60){\special{em:lineto}}
\put(27.43,3.08){\special{em:lineto}}
\put(26.33,3.71){\special{em:lineto}}
\put(25.23,4.55){\special{em:lineto}}
\put(24.38,5.38){\special{em:lineto}}
\put(23.71,6.21){\special{em:lineto}}
\put(23.26,6.81){\special{em:lineto}}
\put(23.26,2.23){\special{em:moveto}}
\put(23.26,6.96){\special{em:lineto}}
\put(23.21,2.16){\special{em:moveto}}
\put(30.03,2.23){\special{em:lineto}}
\put(19.06,2.11){\special{em:moveto}}
\put(18.13,2.48){\special{em:lineto}}
\put(17.26,2.96){\special{em:lineto}}
\put(16.40,3.60){\special{em:lineto}}
\put(15.73,4.23){\special{em:lineto}}
\put(15.18,4.86){\special{em:lineto}}
\put(14.75,5.50){\special{em:lineto}}
\put(14.40,6.13){\special{em:lineto}}
\put(14.15,6.63){\special{em:lineto}}
\put(14.15,6.63){\special{em:moveto}}
\put(13.60,5.53){\special{em:lineto}}
\put(13.05,4.70){\special{em:lineto}}
\put(12.51,4.06){\special{em:lineto}}
\put(11.83,3.43){\special{em:lineto}}
\put(10.93,2.80){\special{em:lineto}}
\put(10.00,2.31){\special{em:lineto}}
\put(9.48,2.11){\special{em:lineto}}
\put(9.48,2.11){\special{em:moveto}}
\put(19.06,2.11){\special{em:lineto}}
\put(14.15,6.63){\special{em:moveto}}
\put(14.15,2.11){\special{em:lineto}}
\end{picture}
}
\caption{Decompositions of bases}
\vspace{1pt}
(see the angles of fundamental triangles under
decompositions).
\end{center}
\label{base}
\end{figure}

\begin{proof}
Let $s$ be a side of $P$ and $b_1$, $b_2$, $b_3$ be the sides
of $F$.
Since $F$ is a Coxeter polyhedron, no angle of $b_i$ is obtuse.
Therefore, any horizontal line is perpendicular to any vertical line
(where the word "line" means an intersection of $s$ with
some mirror). By the same reason the lines are perpendicular to
the edges of $s$.
At most one line may be perpendicular to
the pair of opposite edges of a quadrilateral
in $\H^2$. Thus, there is at most one vertical line and at most
one horizontal line.

\end{proof}
\begin{lemma}\label{l_base}
The bases of $P$ are tiled as shown in Fig.~\ref{base}.

\end{lemma}

\begin{proof}
Let $s$ be a base of $P$.
By Lemma~\ref{krishki} $s$ is either fundamental or a
quasi-Coxeter.  The Coxeter decompositions of triangles are
classified in~\cite{Pink}.
By Lemma~\ref{l_boka} no edge of the triangle $s$ may be
decomposed into more than two parts.  The decompositions
satisfying this condition are listed in Fig.~\ref{base}.


\end{proof}

\begin{theorem}
Let $P$ be a quasi-Coxeter triangular prism such that $F$ is a
triangular prism.  Then the decomposition is one of the
decompositions listed in Table~\ref{prism}.

\end{theorem}

\begin{proof}
It follows from Lemma~\ref{lat}
that two bases of $P$ are decomposed in the same way.
By Lemma~\ref{l_base} there are five possibilities for
decompositions of bases.
Consider two cases.

1). Suppose that the bases are fundamental.
Then the sides are either fundamental or decomposed
like in Fig.~\ref{boka}(c).
Obviously all the sides of $P$ are decomposed in the
same way. This leads to the decomposition shown in the upper part
of Table~\ref{prism}
(since any mirror intersects a boundary of $P$,
there cannot be another decomposition).

2). Suppose that the bases are decomposed. Each decomposition of
the base corresponds to two possible decompositions of the
sides (with the horizontal line or without).
Thus, each decomposition of the base corresponds to
two decompositions of $P$. Some of the dihedral angles of $F$ are
uniquely determined by the combinatorial structure of the
decomposition. We need only to check that the rest dihedral angles
of $F$ can be prescribed in such a way that the dihedral
angles satisfies the condition of Andreev's theorem(see~\cite{A1}).
If the bases are decomposed as in
Fig.~\ref{base}(a) or \ref{base}(c), then it is possible to satisfy
the Andreev's theorem. Otherwise, it is impossible.

\end{proof}

\section{ Decompositions of prisms into tetrahedra}

In this section $P$ is a triangular prism and $F$ is a
tetrahedron.

Any face of $P$ is obviously tiled by triangles,
but in most cases this tiling is not a Coxeter
decomposition.
At first, we derive some information about tilings of bases.

\begin{define}

Let $ABCD$ be a fundamental tetrahedron. An edge $AB$ is called
{\bf $k$-edge} if the dihedral angle formed by $ABC$ and $ABD$
equals $\frac{\pi}{k}$.

\end{define}

\begin{lemma}\label{tile}
Let $p$ be a base of $P$. If any flat angle of $p$ is
fundamental then the tiling of $p$ consists of a unique
triangle. If $p$ has a flat angle decomposed into three parts
then $p$ is bounded by 2-edge, 3-edge and 5-edge.

\end{lemma}

\begin{proof}
Let $f$ be a face of $F$ and $\bar f $ be a plane containing
$f$. Obviously, $\bar f $ is tiled by triangles congruent
to the faces of $F$. To find this tiling consider a face $f$
and the adjacent to $f$ triangles $f_1$, $f_2$ and $f_3$ at the
plane $\bar f$.
These triangles $f_i$ are congruent either to
$f$ or to the other faces of $F$ (we obtain the same face if the
dihedral angle is equal to $\frac{\pi}{k}$, where $k$ is even,
otherwise we obtain another face). Adding face by face we can
prolong the tiling. The only problem is how many triangles are
incident to a fixed vertex.

To solve this problem suppose that $A$ is a vertex of $F$
incident to $k$-edge, $l$-edge and $m$-edge. Suppose that
$\bar f$ contains $A$. Then a decomposition of a small
three-dimensional spherical neighborhood of $A$ is similar to
a Coxeter decomposition of a sphere with fundamental triangle
($\frac{\pi}{k},\frac{\pi}{l},\frac{\pi}{m}$).
An intersection of the neighborhood with $\bar f$ corresponds
to a spherical line in the decomposition of the sphere.
So, to find how many rays starting from $A$ belong to $\bar
f$, it is sufficient to count a number of vertices incident to a
corresponding spherical line.

Thus, for any face of any Coxeter tetrahedron we can find a
tiling of the corresponding plane. In these tilings any triangle
with fundamental angles is a face of $F$. This proves the first
part of the lemma. A triangle with an angle decomposed into
three parts was found in one of these tilings only.
This triangle is bounded by
2-edge, 3-edge and 5-edge and tiled as shown in
Fig.~\ref{treug}.

\begin{figure}[!h]
\begin{center}
{\scriptsize
\setlength{\unitlength}{0.075in}
\begin{picture}(21.46,12.33)
\put(9.06,9.00){{\setbox0=\hbox{3}\lower\ht0\box0}}
\put(19.06,6.33){{\setbox0=\hbox{5}\lower\ht0\box0}}
\put(15.86,3.66){{\setbox0=\hbox{3}\lower\ht0\box0}}
\put(12.26,3.83){{\setbox0=\hbox{2}\lower\ht0\box0}}
\put(9.06,4.50){{\setbox0=\hbox{2}\lower\ht0\box0}}
\put(3.73,5.16){{\setbox0=\hbox{3}\lower\ht0\box0}}
\put(5.20,0.00){{\setbox0=\hbox{2}\lower\ht0\box0}}
\special{em:linewidth 0.014in}
\put(10.66,0.66){\special{em:moveto}}
\put(9.28,2.31){\special{em:lineto}}
\put(7.88,4.51){\special{em:lineto}}
\put(7.20,5.83){\special{em:lineto}}
\put(15.86,12.33){\special{em:moveto}}
\put(12.00,4.50){\special{em:lineto}}
\put(10.66,0.66){\special{em:lineto}}
\put(15.86,12.33){\special{em:moveto}}
\put(15.33,5.50){\special{em:lineto}}
\put(15.46,0.66){\special{em:lineto}}
\put(15.86,12.33){\special{em:moveto}}
\put(9.33,7.33){\special{em:lineto}}
\put(0.00,0.66){\special{em:lineto}}
\put(21.46,0.66){\special{em:moveto}}
\put(19.13,4.36){\special{em:lineto}}
\put(16.86,9.30){\special{em:lineto}}
\put(15.86,12.33){\special{em:lineto}}
\put(0.00,0.66){\special{em:moveto}}
\put(21.46,0.66){\special{em:lineto}}
\end{picture}
}
\end{center}
\caption{The edge labeled by $k$ is a $k$-edge ($k=2,3,5$). }
\label{treug}
\end{figure}
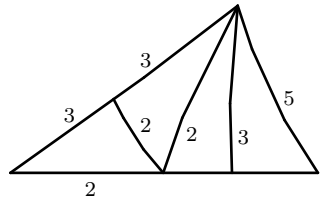

\end{proof}

\begin{lemma}\label{l_imp}
Let $A$ be a non-fundamental vertex of $P$.
Then there exists a non-fundamental edge incident to $A$.

\end{lemma}

\begin{proof}
Suppose that any edge incident to $A$ is fundamental
(i.e. three dihedral angles incident to $A$ are fundamental).
Consider a small sphere $s$ centered in $A$. The decomposition
of $P$ restricted to $s$ is a Coxeter decomposition of a
spherical triangle $p=s\cap P$. Evidently, any angle of $p$
is fundamental.
It is easy to check that in this condition
$p$ is decomposed as shown in Fig.~\ref{fund_s}(a).
Any angle of $p$ is a right angle.

\begin{figure}[!h]
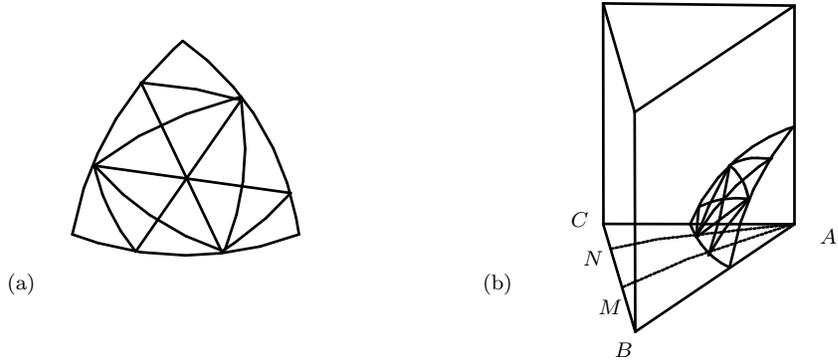

\begin{center}
{\scriptsize
\input pic/imp.tex
}
\caption{If $A$ is a non-fundamental vertex
 but the edges ended in $A$ are fundamental,
then  $AB$ and $AC$ are $2$-edges and the tiling of $\triangle ABC$
is impossible. }
\label{fund_s}
\end{center}
\end{figure}

Thus the neighborhood of $A$ is decomposed as shown in
Fig.~\ref{fund_s}(b). Any edge of $P$ incident to $A$ is a
2-edge. Consider a base $ABC$ of the prism $P$. The angle $A$ of
this triangle is decomposed into three parts by the lines $AN$
and $AM$. It follows from Lemma~\ref{tile} that the sides
of $ABC$ should be 2-edge, 3-edge and 5-edge. This contradicts
to the fact that $AB$ and $AC$ are 2-edges.

\end{proof}

\begin{lemma}
Let $P$ be a triangular prism admitting a Coxeter decomposition
into tetrahedra. Then $P$ has a non-fundamental dihedral angle.

\end{lemma}

\begin{proof}
Suppose that any dihedral angle of $P$ is fundamental.
Then by Lemma~\ref{l_imp}
any
vertex of $P$ is fundamental. By Lemma~\ref{tile} any base of $P$ is
congruent to a single face of $F$. Consider a fundamental tetrahedron
$F_0$ containing a base $\alpha$. Since any dihedral angle of $P$ is
fundamental, the sides of $P$ are the faces of  $F_0$. This is
impossible.

\end{proof}

\subsection*{How to find all decompositions of prisms into tetrahedra }

We have already proved that any prism with tetrahedral
fundamental polyhedron has a non-fundamental dihedral angle.
It follows from Fig.~\ref{mirrors}(1)--(4) that any prism
consists of the tetrahedra, the smaller prisms and
the quadrilateral pyramid. If we know the decompositions of the smaller
parts we can find the decomposition of the whole prism.

Recall that a prism $P$ is called {\bf minimal} if $P$
is non-fundamental and any prism inside $P$ is fundamental.

\begin{define}
We say that a minimal prism is a {\bf prism of level 0}.
A non-fundamental prism $P$ is a {\bf prism of level $k+1$}
if $P$ contains a prism of level $k$ but any prism
inside $P$ contains no prism of level $k$.

\end{define}

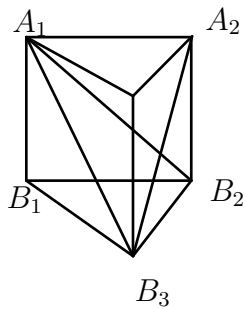
\begin{figure}[!h]
\begin{center}
\setlength{\unitlength}{0.075in}
\begin{picture}(14.15,19.01)
\put(14.15,7.00){{\setbox0=\hbox{$B_2$}\lower\ht0\box0}}
\put(8.95,0.00){{\setbox0=\hbox{$B_3$}\lower\ht0\box0}}
\put(0.01,6.50){{\setbox0=\hbox{$B_1$}\lower\ht0\box0}}
\put(13.88,19.00){{\setbox0=\hbox{$A_2$}\lower\ht0\box0}}
\put(0.41,18.83){{\setbox0=\hbox{$A_1$}\lower\ht0\box0}}
\special{em:linewidth 0.014in}
\put(12.95,17.00){\special{em:moveto}}
\put(8.81,1.66){\special{em:lineto}}
\put(8.81,1.66){\special{em:moveto}}
\put(1.36,16.95){\special{em:lineto}}
\put(12.88,6.93){\special{em:lineto}}
\put(12.88,16.95){\special{em:moveto}}
\put(12.88,6.93){\special{em:lineto}}
\put(8.81,12.83){\special{em:moveto}}
\put(8.81,1.66){\special{em:lineto}}
\put(1.36,16.95){\special{em:moveto}}
\put(1.36,6.93){\special{em:lineto}}
\put(1.36,6.93){\special{em:moveto}}
\put(12.88,6.93){\special{em:lineto}}
\put(8.81,1.66){\special{em:lineto}}
\put(1.36,6.93){\special{em:lineto}}
\put(1.36,16.95){\special{em:moveto}}
\put(12.88,16.95){\special{em:lineto}}
\put(8.81,12.83){\special{em:lineto}}
\put(1.36,16.95){\special{em:lineto}}
\end{picture}
\caption{A minimal prism consists of three (possibly
non-fundamental) tetrahedra.}
\label{3tetr}
\end{center}
\end{figure}

At first, we  classify  decompositions of  minimal
prisms. Evidently, a minimal prism contains no mirrors shown in
Fig.~\ref{mirrors}(2)--\ref{mirrors}(4).
Thus, it contains a mirror shown in Fig.~\ref{mirrors}(1).
Let this mirror be $A_1B_2B_3$ in Fig.~\ref{3tetr}.
Then $A_1$ is a non-fundamental vertex,  and
by Lemma~\ref{l_imp} there is a mirror containing an edge ended in $A_1$.
This edge cannot be $A_1B_1$, since the prism is minimal.
Hence, one of the edges $A_1A_2$ and $A_1A_3$ (say, $A_1A_2$)
is non-fundamental.
Therefore, the prism is decomposed into three tetrahedra.
Decompositions of tetrahedra are classified in~\cite{arxiv-tetr}.
Thus, we can classify decompositions of minimal prisms.

Suppose that we know the classification of decompositions for
the prisms of the levels smaller than $k+1$.  Now we can show how to
find a classification for the prisms of the level $k+1$. Suppose
that $P$ has a dihedral angle decomposed as shown in
Fig.~\ref{mirrors}(2)--\ref{mirrors}(4).  By the definition,
any prismatic part is a prism of the level smaller than $k+1$.
Thus, we know the list of possible decompositions of
prismatic parts. Decompositions of tetrahedral parts
are known too. Therefore, we can classify all
possible decompositions of $P$. Suppose now that $P$ has no
mirrors shown in Fig.~\ref{mirrors}(2)--\ref{mirrors}(4).
Then $P$ has a mirror shown in Fig.~\ref{mirrors}(1).
Moreover, there are  two mirrors of this type decomposing
 $P$ into three tetrahedra
(we use Lemma~\ref{l_imp} again). Combining the decompositions
of these tetrahedra we can find all possible decompositions of
$P$.

Thus, we have an algorithm leading to the classification of
decompositions for the prisms which level is smaller than any
fixed number. In fact, there is no prism of level eight;
hence, there cannot be prism of greater level. Therefore,
this algorithm classifies all decompositions of bounded
triangular prisms. See Table~2 for the classification.

\clearpage

\section*{Tables}

Table~1 contains all decompositions of bounded triangular
prisms which fundamental polyhedron is a triangular prism.
Table~2 contains the list of decompositions of convex bounded
pyramids and triangular prisms which fundamental polyhedron is a
tetrahedron. Here we describe the structure of Table~2.

\begin{itemize}
\item
Horizontal lines separate polyhedra with different
fundamental tetrahedra. For any fundamental tetrahedron the
table contains three columns. The left column contains the list of
decompositions of bounded tetrahedra
(the tetrahedron number 0 is the fundamental one).
The right column lists decompositions of bounded
triangular prisms. The column at the middle contains  bounded
pyramids. At first we list  quadrilateral pyramids which can
be combined from two (possibly non-fundamental) tetrahedra.
Then, after a dotted line, we list the rest quadrilateral
pyramids. After a new dotted line we list pentagonal
pyramids and at last we list hexagonal pyramids.
(There is no Coxeter decomposition of heptagonal pyramid).

\end{itemize}

The polyhedra and their decompositions are represented in
Table~2 as it is described below.

\begin{itemize}

\item
A tetrahedron with  dihedral angles
$\frac{k_i\pi}{q_i}$ $i=1$,...,6 is represented by  the
following diagram:
the nodes of the diagram correspond to the faces of the tetrahedron,
two nodes are connected by a $q$-fold edge decomposed into $k$ parts
if the corresponding faces form up a dihedral angle $\frac{k\pi}{q}$.\\
The numerator $k_i$ is a number of parts in the dihedral angle
(for example, a fraction $\frac{2}{4}$ denotes a right angle
decomposed into two parts).
The faces of the tetrahedron are numbered  by the following way:
the nodes of the diagram have the numbers 0,1,2,3 from the left to the
right.  (The numeration will be used below).

\item
A quadrilateral pyramid $OA_1A_2A_3A_4$ with the base $A_1A_2A_3A_4$
is represented by eight dihedral angles
$$(\widehat{A_1A_2},\widehat{A_2A_3},\widehat{A_3A_4},\widehat{A_4A_1};
\widehat{OA_1},\widehat{OA_2},\widehat{OA_3},\widehat{OA_4}).$$
The dihedral angles are the rational fractions,
where the numerator $q$ means that the corresponding dihedral angle
is decomposed into $q$ parts (we always omit a multiple $\pi$).
A triangular face $OA_iA_{i+1}$
has a number $i$.

\item
By the same way we represent pentagonal and hexagonal
pyramids (we use ten and twelve angles correspondingly).

\item
A triangular prism $A_1A_2A_3B_1B_2B_3$ (where
$A_1A_2A_3$ and $B_1B_2B_3$ are the bases) is represented by nine
dihedral angles
$$(\widehat {A_3B_3},\widehat {A_1B_1},\widehat {A_2B_2};
\widehat {A_1A_2}; \widehat {A_2A_3},\widehat {A_3A_1};
\widehat {B_1B_2},\widehat  {B_2B_3},\widehat {B_3B_1}).$$
The dihedral angles are the rational fractions, where
the numerator is a number of parts in the dihedral
angle and the multiple $\pi$ is omitted. The faces are numbered
in the following way: $A_1A_2A_3$ has number 0,
$B_1B_2B_3$ has number 1,  nd $A_iA_{i+1}B_{i+1}B_i$ has number
$i+1$.

\item
Any decomposition which is a superposition of some other decompositions
is labeled by a star.

\end{itemize}

\Remark
In this notation each pyramid or triangular prism can be written
by several ways (for instance, the order of angles changes
while we rotate the polyhedron). To unify the notation we
choose a record which stands first in the alphabetical order.

\begin{itemize}

\item
For any polyhedron except fundamental tetrahedra we
need some way to reconstruct the decomposition.
For this aim we put some numbers under the diagram of the
polyhedron. These seven numbers (denoted by
$t$,$k$,$l$ and  $m$,$n$,$p$,$q$) show the following:\\
$t$ is a number of the polyhedron (for any fundamental
tetrahedron we start a new numeration),\\
$k$ is a number of fundamental tetrahedra in the decomposition,\\
$l$ is a number of gluings used to obtain this decomposition
(if the decomposition was obtained by a gluing of two polyhedra
with $l=l_1$ and $l=l_2$  then
 $l=1+\max\{l_1,l_2\}$),\\
$m$ and $n$ are the numbers of the polyhedra which should be
glued together to obtain the decomposition.\\
$p$ and $q$ are the numbers of the glued faces of the polyhedra
$m$ and $n$ respectively,\\
The numbers $m$ and $n$ are accompanied by the labels:
we write "tet", "pyr" or "pri"
to show that a corresponding polyhedron
is a tetrahedron, a quadrilateral pyramid or a triangular prism.
(Any tetrahedron admitting a Coxeter decomposition consists of
two smaller tetrahedra,
so we omit the label "tet" in the left column.)

\end{itemize}

\begin{table}[!h]
\caption{Coxeter decompositions of prisms into prisms.}
\begin{center}
\setlength{\unitlength}{0.075in}
\begin{picture}(50.28,99.20)
\put(0.83,45.13){{\setbox0=\hbox{k,m=4
or
5}\lower\ht0\box0}}
\put(13.68,65.65){{\setbox0=\hbox{m}\lower\ht0\box0}}
\put(0.01,66.11){{\setbox0=\hbox{k}\lower\ht0\box0}}
\put(9.06,64.18){{\setbox0=\hbox{2}\lower\ht0\box0}}
\put(12.41,30.23){{\setbox0=\hbox{5}\lower\ht0\box0}}
\put(15.28,33.40){{\setbox0=\hbox{3}\lower\ht0\box0}}
\put(40.00,70.15){{\setbox0=\hbox{3}\lower\ht0\box0}}
\put(41.33,54.75){{\setbox0=\hbox{2}\lower\ht0\box0}}
\put(14.66,21.41){{\setbox0=\hbox{2}\lower\ht0\box0}}
\put(26.83,27.41){{\setbox0=\hbox{4}\lower\ht0\box0}}
\special{em:linewidth 0.014in}
\put(16.98,18.18){\special{em:moveto}}
\put(28.33,17.96){\special{em:lineto}}
\put(25.06,23.56){\special{em:lineto}}
\put(28.16,30.86){\special{em:moveto}}
\put(16.98,30.86){\special{em:lineto}}
\put(16.98,5.28){\special{em:lineto}}
\put(28.33,5.06){\special{em:lineto}}
\put(24.90,10.66){\special{em:lineto}}
\put(24.90,36.23){\special{em:lineto}}
\put(28.16,30.86){\special{em:lineto}}
\put(36.33,61.46){{\setbox0=\hbox{5}\lower\ht0\box0}}
\put(42.88,62.20){{\setbox0=\hbox{2}\lower\ht0\box0}}
\put(44.05,68.21){{\setbox0=\hbox{2}\lower\ht0\box0}}
\put(43.46,53.05){{\setbox0=\hbox{2}\lower\ht0\box0}}
\put(45.23,65.95){\special{em:moveto}}
\put(38.61,65.95){\special{em:lineto}}
\put(38.61,50.35){\special{em:lineto}}
\put(45.23,50.55){\special{em:lineto}}
\put(43.45,53.80){\special{em:lineto}}
\put(43.61,69.16){\special{em:lineto}}
\put(44.43,67.55){\special{em:lineto}}
\put(45.23,65.95){\special{em:lineto}}
\put(40.95,52.80){{\setbox0=\hbox{2}\lower\ht0\box0}}
\put(44.63,59.80){{\setbox0=\hbox{4}\lower\ht0\box0}}
\put(41.15,68.00){{\setbox0=\hbox{2}\lower\ht0\box0}}
\put(20.08,34.61){{\setbox0=\hbox{2}\lower\ht0\box0}}
\put(13.71,23.56){\special{em:moveto}}
\put(28.16,18.18){\special{em:lineto}}
\put(13.61,36.18){\special{em:moveto}}
\put(28.15,30.83){\special{em:lineto}}
\put(28.33,5.28){\special{em:lineto}}
\put(13.91,10.66){\special{em:lineto}}
\put(21.36,21.83){{\setbox0=\hbox{2}\lower\ht0\box0}}
\put(20.08,31.40){{\setbox0=\hbox{2}\lower\ht0\box0}}
\put(21.05,19.03){{\setbox0=\hbox{2}\lower\ht0\box0}}
\put(13.65,23.58){\special{em:moveto}}
\put(36.36,23.58){\special{em:lineto}}
\put(20.30,12.83){\special{em:lineto}}
\put(13.65,23.58){\special{em:lineto}}
\put(36.36,36.33){\special{em:moveto}}
\put(36.51,10.66){\special{em:lineto}}
\put(20.30,25.53){\special{em:moveto}}
\put(20.33,0.23){\special{em:lineto}}
\put(13.65,36.33){\special{em:moveto}}
\put(13.61,10.66){\special{em:lineto}}
\put(13.65,10.73){\special{em:moveto}}
\put(36.36,10.73){\special{em:lineto}}
\put(20.30,0.00){\special{em:lineto}}
\put(13.65,10.73){\special{em:lineto}}
\put(13.65,36.33){\special{em:moveto}}
\put(36.36,36.33){\special{em:lineto}}
\put(20.30,25.53){\special{em:lineto}}
\put(13.65,36.33){\special{em:lineto}}
\put(36.48,69.16){\special{em:moveto}}
\put(45.23,65.95){\special{em:lineto}}
\put(45.23,50.35){\special{em:lineto}}
\put(36.65,53.80){\special{em:lineto}}
\put(50.18,69.25){\special{em:moveto}}
\put(50.26,53.80){\special{em:lineto}}
\put(40.50,62.76){\special{em:moveto}}
\put(40.53,47.51){\special{em:lineto}}
\put(36.50,69.25){\special{em:moveto}}
\put(36.48,53.80){\special{em:lineto}}
\put(36.50,53.85){\special{em:moveto}}
\put(50.18,53.85){\special{em:lineto}}
\put(40.50,47.36){\special{em:lineto}}
\put(36.50,53.85){\special{em:lineto}}
\put(36.50,69.25){\special{em:moveto}}
\put(50.18,69.25){\special{em:lineto}}
\put(40.50,62.76){\special{em:lineto}}
\put(36.50,69.25){\special{em:lineto}}
\put(42.23,83.05){{\setbox0=\hbox{2}\lower\ht0\box0}}
\put(44.80,89.66){{\setbox0=\hbox{2}\lower\ht0\box0}}
\put(42.05,97.20){{\setbox0=\hbox{2}\lower\ht0\box0}}
\put(5.71,68.05){{\setbox0=\hbox{2}\lower\ht0\box0}}
\put(0.70,61.78){\special{em:moveto}}
\put(9.45,58.35){\special{em:lineto}}
\put(0.70,69.48){\special{em:moveto}}
\put(9.45,66.25){\special{em:lineto}}
\put(9.28,50.66){\special{em:lineto}}
\put(0.88,54.10){\special{em:lineto}}
\put(6.16,60.35){{\setbox0=\hbox{2}\lower\ht0\box0}}
\put(10.78,60.56){{\setbox0=\hbox{2}\lower\ht0\box0}}
\put(7.33,62.75){{\setbox0=\hbox{2}\lower\ht0\box0}}
\put(0.73,61.88){\special{em:moveto}}
\put(14.38,61.88){\special{em:lineto}}
\put(4.70,55.41){\special{em:lineto}}
\put(0.73,61.88){\special{em:lineto}}
\put(14.38,69.55){\special{em:moveto}}
\put(14.46,54.10){\special{em:lineto}}
\put(4.70,63.06){\special{em:moveto}}
\put(4.73,47.81){\special{em:lineto}}
\put(0.73,69.55){\special{em:moveto}}
\put(0.70,54.10){\special{em:lineto}}
\put(0.73,54.15){\special{em:moveto}}
\put(14.38,54.15){\special{em:lineto}}
\put(4.70,47.66){\special{em:lineto}}
\put(0.73,54.15){\special{em:lineto}}
\put(0.73,69.55){\special{em:moveto}}
\put(14.38,69.55){\special{em:lineto}}
\put(4.70,63.06){\special{em:lineto}}
\put(0.73,69.55){\special{em:lineto}}
\put(37.23,98.78){\special{em:moveto}}
\put(45.53,95.71){\special{em:lineto}}
\put(45.53,80.95){\special{em:lineto}}
\put(41.46,82.58){\special{em:lineto}}
\put(37.40,84.20){\special{em:lineto}}
\put(50.23,98.85){\special{em:moveto}}
\put(50.28,84.20){\special{em:lineto}}
\put(41.03,92.70){\special{em:moveto}}
\put(41.08,78.25){\special{em:lineto}}
\put(37.26,98.85){\special{em:moveto}}
\put(37.23,84.20){\special{em:lineto}}
\put(37.26,84.25){\special{em:moveto}}
\put(50.23,84.25){\special{em:lineto}}
\put(41.03,78.11){\special{em:lineto}}
\put(37.26,84.25){\special{em:lineto}}
\put(37.26,98.85){\special{em:moveto}}
\put(50.23,98.85){\special{em:lineto}}
\put(41.03,92.70){\special{em:lineto}}
\put(37.26,98.85){\special{em:lineto}}
\put(3.28,89.55){{\setbox0=\hbox{2}\lower\ht0\box0}}
\put(9.48,89.60){{\setbox0=\hbox{2}\lower\ht0\box0}}
\put(8.75,92.55){{\setbox0=\hbox{2}\lower\ht0\box0}}
\put(1.91,91.95){\special{em:moveto}}
\put(14.85,91.95){\special{em:lineto}}
\put(5.70,85.80){\special{em:lineto}}
\put(1.91,91.95){\special{em:lineto}}
\put(14.85,99.20){\special{em:moveto}}
\put(14.88,84.55){\special{em:lineto}}
\put(5.70,93.05){\special{em:moveto}}
\put(5.71,78.61){\special{em:lineto}}
\put(1.91,99.20){\special{em:moveto}}
\put(1.88,84.55){\special{em:lineto}}
\put(1.91,84.63){\special{em:moveto}}
\put(14.85,84.63){\special{em:lineto}}
\put(5.70,78.45){\special{em:lineto}}
\put(1.91,84.63){\special{em:lineto}}
\put(1.91,99.20){\special{em:moveto}}
\put(14.85,99.20){\special{em:lineto}}
\put(5.70,93.05){\special{em:lineto}}
\put(1.91,99.20){\special{em:lineto}}
\end{picture}
\end{center}
\label{prism}
\end{table}

\clearpage

\begin{center}
{\normalsize
Table~2. Coxeter decompositions of  pyramids and prisms\\
which fundamental polyhedron is a tetrahedron}

\vspace{15pt}

\input pic/tabpr4e.tex
\end{center}


\pagebreak

\end{document}